\documentclass[12pt]{amsart}
\usepackage{amsfonts,amssymb}
\usepackage{graphicx}
\usepackage{amscd}
\input epsf

\newtheorem{thm}{Theorem}[section]
\newtheorem{prop}[thm]{Proposition}
\newtheorem{cor}[thm]{Corollary}
\newtheorem{lemma}[thm]{Lemma}
\newtheorem{conj}[thm]{Conjecture}

\theoremstyle{definition}

\theoremstyle{remark}
\newtheorem{remark}{Remark}

\def\proof{{\sc Proof. }}
\def\A{\mathbb{A}}
\def\R{\mathbb{R}}
\def\C{\mathbb{C}}
\def\Z{\mathbb{Z}}
\def\H{\mathbb{H}}

\def\Lc{{\mathcal L}}
\def\tr{{\rm tr}}
\def\det{{\rm det}}
\def\Mat{{\rm Mat}}

\def\SL{{\rm SL}}
\def\<{\langle }
\def\>{\rangle }
\def\c{{\rm c}}
\def\s{{\rm s}}
\def\eps{\epsilon}
\def\Box{\square}

\def\x{\mathbf{x}}
\def\y{\mathbf{y}}
\def\z{\mathbf{z}}
\def\e{\mathbf{e}}
\def\w{\mathbf{w}}

\title{On binary quadratic forms with semigroup property}

\author[F. Aicardi]{Francesca Aicardi}
\email{<aicardi@sissa.it>}

\author[V. Timorin]{Vladlen Timorin}
\email{<timorin@math.sunysb.edu>}

\begin{document}

\begin{abstract}
A quadratic form $f$ is said to have {\em semigroup property} if its values at
points of the integer lattice form a semigroup under multiplication.
A problem of V. Arnold is to describe all binary integer
quadratic forms with semigroup property.
If there is an integer bilinear map $s$ such that
$f(s(\x,\y))=f(\x)f(\y)$ for all vectors $\x$ and $\y$ from the integer
2-dimensional lattice, then the form $f$ has semigroup property.
We give an explicit integer parameterization of all pairs $(f,s)$ with the
property stated above.
We do not know any other examples of forms with semigroup property.
\end{abstract}

\maketitle

\section*{Introduction}

We are dealing with {\em binary quadratic forms}, i.e. real-valued quadratic forms
defined on $\R^2$ (2-dimensional real coordinate space).
From now on, the word ``binary'' will be suppressed.
An {\em integer quadratic form} is a quadratic form on $\R^2$ that takes
integer values at all points of the integer lattice $\Z^2\subset\R^2$.
In coordinates, an integer quadratic form is given by a formula
$mx_1^2+kx_1x_2+nx_2^2$, where the coefficients $m$, $k$ and $n$ are integers.
This form is abbreviated as $(m,k,n)$.
An integer quadratic form is said to have {\em semigroup property} if the
product of its values at any two integer points is also a value at an
integer point\footnote{In \cite{Ai,Ar}, integer quadratic forms with semigroup property
are called {\em perfect forms}.
However, the term ``perfect quadratic form/Euclidean lattice'' already existed in
the literature (it dates back to works of G. Voronoi written in early 1900s, for
a modern account see \cite{Mart}).
The authors are grateful to P. Gunnells for pointing this out.}.
In other words, a quadratic form $f:\R^2\to\R$ has semigroup property if for any
two vectors $\x,\y\in\Z^2$ there exists a third vector $\z\in\Z^2$ such that $f(\z)=f(\x)f(\y)$.

Let us start by giving two very simple examples that were known to
ancient mathematicians.

\subsection*{Example}
Consider the quadratic form $(1,0,1)$.
This form has semigroup property, as we can see from the following formula:
$$
(x_1^2+x_2^2)(y_1^2+y_2^2)=(x_1y_1-x_2y_2)^2+(x_1y_2+x_2y_1)^2.
$$
This formula has a nice interpretation in term of multiplication of complex numbers.
Namely, identify the plane $\R^2$ with the field $\C$ of complex numbers in the following way.
A vector with coordinates $(x_1,x_2)$ gets identified with the complex number $x_1+ix_2$.
The formula displayed above means exactly that $|x|^2|y|^2=|xy|^2$ for any two
complex numbers $x$ and $y$.

More generally, the form $(1,0,D)$ has semigroup property.
Namely, we have the following formula:
$$
(x_1^2+Dx_2^2)(y_1^2+Dy_2^2)=(x_1y_1-Dx_2y_2)^2+D(x_1y_2+x_2y_1)^2.
$$
For $D>0$, it also has a nice interpretation in terms of complex numbers.
Namely, we should identify a vector $(x_1,x_2)$ with the complex
number $x_1+\sqrt{-D}x_2$.
Then the formula reduces to the identity $|x|^2|y|^2=|xy|^2$ for complex
numbers $x$ and $y$.

Consider a bilinear map $s:\R^2\times\R^2\to\R^2$.
We say that $s$ is an {\em integer normed pairing} with respect to a
quadratic form $f:\R^2\to\R$ if it takes $\Z^2\times\Z^2$ to $\Z^2$
and satisfies the property
$$
f(s(\x,\y))=f(\x)f(\y)
$$
for all $\x,\y\in\R^2$.
Note that we do not require that the form $f$ be an integer form.
As we will see, any non-degenerate quadratic form admitting an integer normed
pairing must be an integer quadratic form.
An integer quadratic form admitting an integer normed pairing, has semigroup
property.
In the examples of forms with semigroup property considered above,
we always had an integer normed pairing.

\begin{conj}
If an integer quadratic form has semigroup property,
then it admits an integer normed pairing.
\end{conj}

V. Arnold \cite{Ar} posed the following problem: describe all forms
with semigroup property.
In this article, we solve a much easier problem: we will describe all
nondegenerate quadratic forms admitting integer normed pairings.

\subsection*{A brief tour}
The main results of this paper, which will be stated later in full detail,
are roughly as follows.
We give explicit polynomial formulas for all possible
integer normed pairings $s:\R^2\times\R^2\to\R^2$.
They split into four families according to their {\em type}.
For a bilinear product $s:\R^2\times\R^2\to\R^2$ and a fixed vector $\y$,
consider the linear map $\x\mapsto s(\x,\y)$.
Divide the determinant of this map by $f(\y)$.
It is not hard to prove that the sign of this ratio is independent of $\y$
whenever the ratio is defined.
Denote this sign by $\epsilon_1$.
Similarly, define $\epsilon_2$ as the sign of $\det(\y\mapsto s(\x,\y))/f(\x)$.
We will say that $s$ is of type $(\epsilon_1;\epsilon_2)=(\pm;\pm)$.

Integer normed pairings of types $(+,+)$, $(-,+)$ and $(+,-)$
depend on 5 integer parameters and always come together:
if a nondegenerate quadratic form admits a normed pairing of one of
these types, then it admits normed pairing of two other types as
well.
However, integer normed pairings of type $(-,-)$, which depend
on 4 integer parameters, behave very differently.

The arithmetic of positive definite quadratic forms is
closely related to geometry of two-dimensional lattices in $\C$.
Namely, quadratic forms considered up to $SL_2(\Z)$-equivalence
correspond to lattices up to rotations.
Multiplication of complex numbers and the quadratic form
$z\mapsto |z|^2$ play an important role here.
We will describe integer normed pairings of types $(+,+)$, $(-,+)$ and $(+,-)$
as those corresponding to ideals in the rings of quadratic integers
(i.e. the corresponding lattices will literally coincide with ideals).
Quadratic forms admitting normed pairings of only one type, $(-,-)$,
give rise to elements of order 3 in the class groups.

Indefinite forms correspond to lattices in the algebra $\H$ of
hyperbolic numbers, i.e. the algebra generated over reals by 1 and
an element $j$ satisfying the equation $j^2=1$.
Indefinite forms with semigroup property admit analogous results
to those for positive definite forms.
We only need to replace $\C$ with $\H$.

Another interesting aspect of integer normed pairings is their
connections with the trigroup property.
V. Arnold pointed out \cite{Ar} that, while not all integer quadratic forms
on $\R^2$ have semigroup property, all of them have {\em trigroup property}:
the product of three integers represented by a given integer quadratic
form, is also represented by this form.
Although a formal corollary from the Gauss composition,
this is a particularly nice and simple statement, worth mentioning
in elementary textbooks.

From the trigroup property it follows, in particular, that
an integer quadratic form representing 1 (i.e., attaining the value 1
somewhere on $\Z^2$) satisfies the semigroup property.
It turns out that such a form admits integer normed pairings of all
4 types in this case.

A more general corollary from the trigroup property is that if
an integer quadratic form $g$ represents an integer $m$, then the
form $f=mg$ has semigroup property \cite{Ar}.
This form always admits integer normed pairings of
types $(+,+)$, $(-,+)$ and $(+,-)$, and sometimes admits a
normed pairing of the type $(-,-)$, but not always.

We will give explicit examples in Section 6.

\subsection*{Acknowledgements.}
The authors are grateful to Vladimir Arnold for generating the interest
in the subject and for useful discussions, and to John Milnor for suggesting
improvements of the exposition.

\section{Statement of main results}

Let us first state the results, and then explain their geometric meaning.

A bilinear pairing $s:\R^2\times\R^2\to\R^2$ is represented in coordinates by a pair of
bilinear forms.
Using matrix product notation, a bilinear product $s(\x,\y)=\z$
can be described by a pair of $2\times 2$ matrices $A_j$ where
$$
z_j=(x_1,x_2)A_j\left(\begin{array}{c}y_1\\y_2\end{array}\right),\quad j=1,2,\
\x=(x_1,x_2),\ \y=(y_1,y_2),\ \z=(z_1,z_2)
$$
We will write $s$ as $(A_1|A_2)$.

\begin{thm}
\label{inps}
Let $f$ be a nondegenerate quadratic form on $\R^2$.
Suppose that there exists an integer normed pairing $s$ with respect to $f$.
Then the pair $(s,f)$ is defined by one of the following formulas:
$$
s_1=\left(\begin{array}{cc}
mp+kq&nq\\
nq&-np
\end{array}\vline
\begin{array}{cc}
-mq&mp\\
mp&nq+kp
\end{array}
\right),\quad
\begin{array}{c}
f_1=(rm,rk,rn),\\
\hbox{where }r:=mp^2+kpq+nq^2.
\end{array}
$$
$$
s_2=\left(\begin{array}{cc}
mp&nq+kp\\
-nq&np
\end{array}\vline
\begin{array}{cc}
mq&-mp\\
mp+kq&nq
\end{array}
\right),\quad
\begin{array}{c}
f_2=(rm,rk,rn),\\
\hbox{where }r:=mp^2+kpq+nq^2.
\end{array}
$$
$$
s_3=\left(\begin{array}{cc}
mp&-nq\\
nq+kp&np
\end{array}\vline
\begin{array}{cc}
mq&mp+kq\\
-mp&nq
\end{array}
\right),\quad
\begin{array}{c}
f_3=(rm,rk,rn),\\
\hbox{where }r:=mp^2+kpq+nq^2.
\end{array}
$$
$$
s_4=\left(\begin{array}{cc}
a & c\\c  & b
\end{array}\vline
\begin{array}{cc}
-d & -a \\ -a & -c
\end{array}\right),\quad
f_4=(a^2-cd,\ ac-bd,\ c^2-ab)
$$
for some integers $m$, $n$, $k$, $p$, $q$,
$a$, $b$, $c$, $d$.
In particular, the form $f$ is an integer form.

Conversely, for arbitrary integers $m$, $n$, $k$, $p$,
$q$, $a$, $b$, $c$, $d$, the
forms $f_1$, $f_2$, $f_3$ and $f_4$ given by the formulas
displayed above admit integer normed pairings
$s_1$, $s_2$, $s_3$ and $s_4$, respectively.
They are of types $(+,+)$, $(-,+)$, $(+,-)$ and $(-,-)$,
respectively.
\end{thm}

%A proof of this theorem is given on page \pageref{p:inps}.
One can see that the first three formulas
(for $f_i$ and $s_i$, where $i=1,2,3$) are similar to each other.
Their integer normed pairings are defined by the same five integer parameters.
The formulas provided for the quadratic forms $f_i$, $i=1,2,3$ are
equivalent to saying that each $f_i$ is of type $rg$, where
an integer quadratic form $g$ represents an integer $r$.
In particular, if $rg$ is primitive, then $r=1$, and $g$
represents the identity element in the class group.

However, the last formula (for $s_4$ and $f_4$) is very different from the first three.
In this context, the formula for $f_4$ and $s_4$ appeared in \cite{Ai}.
It may be regarded as a parameterization of all quadratic forms
having order 1 or 3 in the class group.
As such, it is contained in the remarkable article \cite{Bhar} of
M. Bhargava, who associates the quadratic forms of type $f_4$ with
classes of cubic forms.
The statement that all pairs $(f_4,s_4)$, where $s_4$ is an
integer normed pairing of type $(-,-)$ associated with a quadratic form $f_4$,
are parameterized as above, can be easily proved by comparing
Bhargava's approach to composition with Gauss' original approach.

\subsection*{Geometric meaning of Theorem \ref{inps}.}
There is a remarkable correspondence between positive definite
integer quad\-ra\-tic forms and certain 2-di\-men\-si\-onal lattices in $\C$.
Given a positive definite quadratic form $f$ on $\R^2$, consider an
orientation preserving real vector space isomorphism $\phi:\R^2\to\C$
such that $f(\x)=|\phi(\x)|^2$ for all $\x\in\R^2$.
The orientation in $\R^2$ is given by the standard basis $(1,0)$, $(0,1)$.
The orientation in $\C$ is given by the basis $(1,i)$.
Such isomorphism $\phi$ exists because any positive definite quadratic
form reduces to the sum of squares.

The image $L$ of the integer lattice $\Z^2$ under the isomorphism $\phi$
is a 2-dimensional lattice in $\C$.
This lattice is {\em integer-normed}, i.e. the number $|z|^2$ is an integer
for any $z\in L$.
The theory of positive definite integer quadratic forms is parallel to
the theory of integer-normed lattices in $\C$.
We will say that the lattice $L$ {\em corresponds} to the quadratic form $f$.

We will give a proof of the following theorem
(which is very classical and probably goes back to Gauss;
Hurwitz proved a much more general fact):

\begin{thm}
\label{sigma1-4}
If a positive definite integer quadratic form $f$ on $\R^2$ admits a
normed pairing, then there exists a corresponding lattice $L$ in $\C$
that is stable under at least one of the following operations:
$$
\begin{array}{c}
\sigma_1:(z,w)\mapsto zw,\\
\sigma_2:(z,w)\mapsto \overline{z}w,\\
\sigma_3:(z,w)\mapsto z\overline{w},\\
\sigma_4:(z,w)\mapsto \overline{zw}.
\end{array}
$$
\end{thm}

A proof of this theorem is given on page \pageref{p:sigma1-4}.

Suppose that a lattice $L$ is stable under $\sigma_k$, where $k=1,\dots, 4$.
Consider an orientation preserving vector space isomorphism $\phi:\R^2\to\C$
such that $\phi(\Z^2)=L$.
Then the quadratic form $f=|\phi|^2$ admits the following integer normed pairing:
$$
s(\x,\y)=\phi^{-1}\sigma_k(\phi(\x),\phi(\y)),\quad \x,\y\in\R^2.
$$
We will say that the pairing $s$ {\em corresponds} to the operation $\sigma_k$.

\begin{thm}
All integer normed pairings of types $(+,+)$, $(-,+)$, $(+,-)$ and $(-,-)$
with respect to the quadratic form $f$
correspond to the operations $\sigma_1$, $\sigma_2$, $\sigma_3$ and $\sigma_4$,
respectively.
\end{thm}

A proof of this theorem is given on page \pageref{p:inps},
together with a proof of Theorem \ref{inps}.

The lattices stable under operations $\sigma_1$ through $\sigma_3$ admit the following
simple description.

Suppose that $\Delta$ is a negative number congruent to 0 or 1 modulo 4.
Define the ring $R_\Delta$ as follows.
If $\Delta$ is divisible by 4, then $R_\Delta$ is $\Z[\sqrt{\Delta}/2]$.
If $\Delta$ is congruent to 1 modulo 4, then $R_\Delta$ is $\Z[(1+\sqrt{\Delta})/2]$.
Any ring of quadratic integers (algebraic integers of degree 2) in $\C$ that
contains $\Z$ and is not contained in $\R$,
coincides with $R_\Delta$ for some $\Delta<0$ congruent to 0 or 1 modulo 4.

\begin{thm}
\label{sigma1-3}
Suppose that a lattice in $\C$ is stable under one of the operations
$\sigma_1$, $\sigma_2$, $\sigma_3$.
Then it is stable under all these operations.
Moreover, this lattice coincides with an ideal in the ring $R_\Delta$ for
some negative $\Delta$ congruent to 0 or 1 modulo 4.
In particular, this lattice is integer-normed.
\end{thm}

For a proof of this theorem, see page \pageref{p:sigma1-3}.

This theorem gives the following immediate corollary in terms of
integer quadratic forms:

\begin{cor}
\label{s1-3}
If a positive-definite integer quadratic form admits a normed pairing
of type $(+,+)$, $(-,+)$, or $(+,-)$, then it admits normed pairings of
all three types.
\end{cor}

It is harder to describe lattices stable under $\sigma_4$.
But they have a remarkable property, which is best formulated in terms of the class group.

For two additive subgroups $L'$ and $L''$ in $\C$ (in particular, for two
lattices), we can define the product
$L'L''$ as the additive subgroup of $\C$ generated by all pair-wise products $z'z''$,
where $z'\in L'$ and $z''\in L''$.
In general, the product $L'L''$ thus defined would not be a 2-dimensional lattice:
it can be an everywhere dense subgroup.
However, it turns out that the product of two integer-normed lattices
of the same discriminant is also a lattice.

The {\em discriminant} of a lattice $L$ is defined as the discriminant of
a corresponding quadratic form.
Geometrically, it is equal to the square of the area of a fundamental parallelogram
multiplied by $-4$.

An integer-normed lattice $L\subset\C$ is said to be {\em primitive} if the greatest
common divisor of all numbers $|z|^2$, where $z\in L$, is equal to 1.
It is a very classical fact that the product of two integer-normed lattices
of the same discriminant $\Delta$ is also an integer lattice.
Moreover, if two integer-normed lattices of discriminant $\Delta$ are primitive, then
their product also has discriminant $\Delta$.

Primitive integer-normed lattices of discriminant $\Delta$ form a group under
multiplication.
Let us denote this group by $\Lc(\Delta)$.
The identity element in the group $\Lc(\Delta)$ is exactly the ring $R_\Delta$,
which itself is an integer-normed lattice of discriminant $\Delta$.
The inverse element to a lattice $L\in\Lc(\Delta)$ is the complex conjugate lattice
$$
\overline{L}=\{\overline{z}|\ z\in L\}.
$$

Denote by $S^1$ the multiplicative group of complex numbers of length 1.
Thus $S^1=\{\alpha\in\C\ |\ |\alpha|=1\}$.
There is a natural homomorphism of the group $S^1$ to the group $\Lc(\Delta)$.
Under this homomorphism, a number $\alpha\in S^1$ gets mapped to $\alpha R_\Delta$.
The quotient group of $\Lc(\Delta)$ by the image of $S^1$ is called
{\em the class group}.
It is isomorphic to the ideal class group of the ring $R_\Delta$.
The class group was first introduces by Gauss \cite{Gauss} in terms
of integer quadratic forms.

\begin{thm}
\label{ord3}
Suppose that a primitive integer-normed lattice $L$ is stable under $\sigma_4$.
Then $L^3=R_\Delta$.
If $L\ne R_\Delta$, then $L$ represents an element of order 3 in the group $\Lc(\Delta)$.
If $L$ does not belong to the image of $S^1$ in $\Lc(\Delta)$,
then $L$ represents an element of order 3 in the class group.
\end{thm}

This theorem is almost obvious, so we give a proof here:

\proof
If $L$ is stable under $\sigma_4$, then $\overline{L}\overline{L}=L$.
But $\overline{L}=L^{-1}$ in the group $\Lc(\Delta)$.
Therefore, $L^3=1$ in this group.
$\Box$

Since we have an explicit formula for all integer normed
pairings of type $(-,-)$ (see Theorem \ref{inps}),
Theorem \ref{ord3} might have interesting applications to description of
order 3 elements in class groups.
However, the main difficulty is to determine if a quadratic
form $f_4$ from Theorem \ref{inps} is not principal, i.e. if it
does not correspond to the identity element in the class group.
Some parameterization of all order 3 elements in class groups
is known \cite{Kishi-Miyake2000}, but it is given in different terms.

Suppose now that an integer quadratic form is indefinite
(i.e. it can take positive values as well as negative values,
in particular, it is non-degenerate).
Then, analogously to what we did for positive definite forms, we can
define a corresponding integer lattice in the algebra $\H$ of hyperbolic
numbers.
Recall that the algebra $\H$ of hyperbolic numbers is defined
as the quotient of the polynomial algebra $\R[t]$ by the
principal ideal $(t^2-1)$.
The algebra $\H$ is a 2-dimensional real vector space spanned by 1 and
the image $j$ of the polynomial $t$.
We know that $j$ satisfies the equation $j^2=1$. The algebra of hyperbolic
numbers contains zero divisors, e.g. $(1+j)(1-j)=0$.
By definition, $\H$ is associative and commutative.

There are only two automorphisms of $\H$ over real numbers,
namely, the identity and the conjugation.
The conjugation can be defined as the $\R$-linear map from $\H$ to $\H$ such that
$j$ gets mapped to $-j$.
For a hyperbolic number $z$, denote its conjugate by $\overline{z}$.

If $z=\xi+\eta j\in\H$, where $\xi$ and $\eta$ are real numbers,
then define $|z|^2$ to be $\xi^2-\eta^2=z\overline{z}$.
Thus the norm $|z|$ of $z$ is not always a real number.
It can be a purely imaginary number.

It turns out that all our statements concerning positive definite
quadratic forms and lattices in $\C$ generalize to indefinite
quadratic forms and lattices in $\H$.
However, in order to extend Theorem \ref{sigma1-3}, we need to give
another definition of the ring $R_\Delta$ for positive $\Delta$
congruent to 0 or 1 modulo 4.
If $\Delta$ is divisible by $4$, then define $R_\Delta$ as $\Z[j\sqrt{\Delta}/2]$.
If $\Delta$ is congruent to 1 modulo 4, then define
$R_\Delta$ as $\Z[(1+j\sqrt{\Delta})/2]$.
Thus $R_\Delta$ is not a subring of $\R$.

\begin{thm}
\label{hyper}
Theorems \ref{sigma1-4} generalizes to indefinite forms.
We only need to replace the field $\C$ of complex numbers with the
algebra $\H$ of hyperbolic numbers.
Accordingly, Theorem \ref{sigma1-3} is also true for lattices in $\H$,
and Corollary \ref{s1-3} is also true for indefinite forms.
\end{thm}

We will prove this theorem simultaneously with the corresponding theorems
for positive definite forms.
To this end, we will introduce a unified notation that works equally well
in both cases (positive definite and indefinite).

We have discussed the result of Theorem \ref{inps} in terms of lattices
and ideals.
Let us now describe some interpretations of the pairings $s_i$ in terms
of integer quadratic forms.

\subsection*{Embeddings into the matrix algebra.}
Consider the lattice $\Mat_2(\Z)$ of all integer $2\times 2$-matrices.
This is a 4-dimensional lattice.
This lattice is equipped with a natural quadratic form, namely, the determinant.
Multiplication of matrices is obviously a normed pairing with respect to the determinant.

But there are other normed pairings.
For a matrix $A\in\Mat_2(\Z)$, denote by $\overline{A}$ its adjoint matrix.
Recall that the adjoint matrix $\overline{A}$ is also an integer matrix, and
it satisfies the relation
$A\overline{A}=\overline{A}A=\det(A)E$, where $E$ denotes the identity matrix.
Define the following integer pairings of the lattice $\Mat_2(\Z)$:
\begin{equation}
\label{Si}
\begin{array}{c}
S_1: (A,B)\mapsto AB,\\
S_2: (A,B)\mapsto \overline{A}B,\\
S_3: (A,B)\mapsto A\overline{B},\\
S_4: (A,B)\mapsto \overline{AB},\\
\end{array}
\end{equation}
The pairings $S_1$, $S_2$, $S_3$ and $S_4$ are normed pairings with
respect to the determinant.

To obtain 2-dimensional examples, it suffices
to find 2-dimensional sublattices of the lattice $\Mat_2(\Z)$ that
are stable under one of the pairings $S_1$, $S_2$, $S_3$ or $S_4$.

\begin{thm}
\label{S1-3}
Let $A\in\Mat_2(\Z)$ be any integer $2\times 2$ matrix.
Take an integer $r$ that divides the determinant of $A$.
Then the sublattice of $\Mat_2(\Z)$ spanned by $A$ and $rE$, is stable under the
pairings $S_1$, $S_2$ and $S_3$.
The pairings $S_1$, $S_2$ and $S_3$ restricted to the sublattice thus
constructed are integer normed pairings of types $(+,+)$, $(-,+)$
and $(+,-)$, respectively, with respect to the restriction of the determinant.
Any two-dimensional sublattice of $\Mat_2(\Z)$ not belonging to the zero level
of the determinant and stable under one of the pairings $S_k$, $k=1,2,3$,
can be obtained as above.
\end{thm}

This theorem will be proved in Section 2 (see Propositions \ref{constrS1-3},
\ref{sublatS1}, \ref{sublatS2}, \ref{sublatS3} and Theorem \ref{S1-32s1-3}).
The converse is also true:

\begin{thm}
\label{s1-3toS1-3}
Let $s$ be an integer normed pairing of type $(+,+)$, $(-,+)$ or $(-,-)$.
Then there is a linear embedding $\phi:\R^2\to\Mat_2(\R)$ such that
$\phi(\Z^2)$ is a sublattice of $\Mat_2(\Z)$ and
$$
\phi(s(\x,\y))=S_k(\phi(\x),\phi(\y)),\quad f(\x)=\det(\phi(\x))
$$
for all $\x,\y\in\R^2$, where $k$ is 1, 2 or 3.
\end{thm}

This theorem will follow from Theorem \ref{sigma1-4} and
Theorem \ref{s1-32S1-3} proved in Section 2.

Analogously, sublattices of $\Mat_2(\Z)$ stable under $S_4$ correspond
to normed pairings of type $(-,-)$.

\begin{thm}
\label{S4}
Let $A\in\Mat_2(\Z)$ be any integer $2\times 2$ matrix.
Take an integer $r$ that divides the number $\tr(A)^2-\det(A)$.
Then the sublattice of $\Mat_2(\Z)$ spanned by $A$ and $rE$,
is stable under the pairing $S_4$.
The pairing $S_4$ restricted to the sublattice thus constructed
is an integer normed pairing of type $(-,-)$.
Any 2-dimensional sublattice of $\Mat_2(\Z)$ stable under $S_4$
can be obtained in this way.
\end{thm}

However, we do not obtain all integer pairings of type $(-,-)$
by embedding the integer lattice $\Z^2$ into $\Mat_2(\Z)$ as
a sublattice stable under $S_4$.

Theorem \ref{S4} will be proved in Section 2
(see Propositions \ref{constrS4}, \ref{sublatS4} and Theorem \ref{S42s4}).

\subsection*{Commutative traceless pairings.}
An interpretation of integer normed pairings of type $(-,-)$ can be obtained
as follows.
Consider an integer pairing $s:\R^2\times\R^2\to\R^2$
(i.e. a bilinear map such that $s(\Z^2\times\Z^2)\subseteq\Z^2$)
satisfying the following properties:
\begin{itemize}
\item it is commutative, i.e. $s(\x,\y)=s(\y,\x)$
for all $\x,\y\in\R^2$,
\item it is traceless, i.e. for any $\x\in\R^2$ the
operator $M_\x:\y\mapsto s(\x,\y)$ has trace zero.
\end{itemize}
It turns out that any such pairing is a normed pairing of type $(-,-)$
with respect to some integer quadratic form.

Denote by $\e_1=(1,0)$, $\e_2=(0,1)$ the standard basis in $\R^2$.
Any pairing $s:\R^2\times\R^2\to\R^2$ is determined by the pair of
matrices $M_1=M_{\e_1}$ and $M_2=M_{\e_2}$.
Indeed, for all $\x,\y\in\R^2$ we have
$$
s(\x,\y)=x_1M_1(\y)+x_2M_2(\y).
$$

The conditions of being commutative and traceless are linear.
That is why we can give an explicit formula for all commutative
traceless pairings.
We have
\begin{equation}
\label{cotrl}
M_1=\left(\begin{array}{cc}
a& c\\
-d& -a\end{array}\right),\quad
M_2=\left(\begin{array}{cc}
c& b\\
-a& -c\end{array}\right),
\end{equation}
where $a$,$b$, $c$ and $d$ are arbitrary integers.
The fact that $s$ is traceless means that both matrices
$M_1$ and $M_2$ are traceless.
The commutativity of $s$ means that the second column
of $M_1$ coincides with the first column of $M_2$.
Therefore, the pairing $s$ is commutative and traceless if and only if the
matrices $M_1$ and $M_2$ are given by equation (\ref{cotrl}).

\begin{thm}
Any commutative traceless pairing $s:\R^2\times\R^2\to\R^2$ is a
normed pairing of type $(-,-)$ with respect to a certain quadratic
form $f$, which can be recovered from the relation
$$
s(\x,s(\x,\y))=f(\x)\y,\quad \x,\y\in\R^2
$$.
Any integer normed pairing of type $(-,-)$ is commutative and traceless.
In particular, it is given by formula (\ref{cotrl}).
\end{thm}

\proof
It is easy to verify by direct computation that this description of
$s$ and $f$ gives exactly the formula for $(f_4,s_4)$ from
Theorem \ref{inps}.
The parameters $a$, $b$, $c$, $d$ from formula (\ref{cotrl}) match the
parameters $a$, $b$, $c$, $d$ form Theorem \ref{inps}.
$\Box$

\subsection*{Plan.}
This article is organized as follows.
In Section 2, we describe sublattices of $\Mat_2(\Z)$.
In particular, we prove Theorems \ref{S1-3} and \ref{S4}.
In Section 3, we recall the trigroup property discovered by V. Arnold,
and explain relationships between this property and integer normed pairings.
We will give a new proof of the trigroup property.
Sections 4 and 5 deal with lattices in $\C$ or in $\H$ stable under
operations $\sigma_k$, $k=1,\dots, 4$.
We will prove all main theorems there.
Finally, Section 6 contains some examples concerning integer normed pairings.

\section{Sublattices in the matrix algebra}

An element $x$ of an abstract lattice $\Lambda$ is called {\em primitive}
if it does not have the form $ay$ for an element $y\in\Lambda$ and
an integer $a>1$.

Recall the following well-known fact about 2-dimensional lattices.

\begin{prop}
\label{basis}
Any primitive element of a 2-dimensional lattice is contained in
some basis of the lattice.
\end{prop}

\proof Identify the lattice with $\Z^2$. A vector $\x\in\Z^2$ with
coordinates $(x_1,x_2)$ is primitive if and only if $x_1$ and
$x_2$ are relatively prime. Then there are integers $y_1$ and $y_2$
such that $x_1y_2-x_2y_1=1$. This means that the vectors
$(x_1,x_2)$ and $(y_1,y_2)$ form a basis of $\Z^2$. $\Box$

Recall that the pairings $S_1$ through $S_4$ on $\Mat_2(\Z)$ are normed
pairings with respect to the determinant.

Let us say that a sublattice $\Lambda$ of $\Mat_2(\Z)$ is a {\em null sublattice}
if the determinants of all matrices from $\Lambda$ are zero.
Otherwise a sublattice $\Lambda$ is said to be a {\em non-null sublattice}.

We will first explain how to construct sublattices of $\Mat_2(\Z)$
stable under the pairings $S_k$, $k=1,\dots, 4$, and then prove that our
construction gives all such non-null sublattices.

\begin{prop}
\label{constrS1-3}
Let $A\in\Mat_2(\Z)$ be any integer $2\times 2$ matrix.
Take an integer $r$ that divides the determinant of $A$.
Then the sublattice of $\Mat_2(\Z)$ spanned by $A$ and $rE$, is
stable under the pairings $S_1$, $S_2$ and $S_3$.
\end{prop}

\proof Denote by $\Lambda$ the lattice spanned by $A$ and $rE$. It
suffices to show that the matrices $A^2$, $r\overline{A}$ and
$A\overline{A}$ belong to the lattice $\Lambda$.

By the Cayley--Hamilton theorem, we have
$$
A^2-\tr(A)A+\det(A)E=0.
$$
Therefore, the matrix $A^2$ belongs to the lattice $\Lambda$. But
we also have
$$
A+\overline{A}=\tr(A)E.
$$
This means that $r\overline{A}$ also belongs to $\Lambda$. The
product $A\overline{A}=\det(A)E$ is in $\Lambda$, since the
determinant of $A$ divisible by $r$. This concludes the proof.
$\Box$

\begin{prop}
\label{constrS4}
Let $A\in\Mat_2(\Z)$ be any integer $2\times 2$ matrix.
Take an integer $r$ that divides the number $\tr(A)^2-\det(A)$.
Then the sublattice of $\Mat_2(\Z)$ spanned by $A$ and $rE$,
is stable under the pairing $S_4$.
\end{prop}

\proof
It suffices to verify that the matrices $\overline{A}^2$
and $r\overline{A}$ belong to the sublattice $\Lambda$ spanned by
$A$ and $rE$.
Indeed,
$$
\begin{array}{c}
\overline{A}^2=\tr(A)\overline{A}-\det(A)E=-\tr(A)A+(\tr(A)^2-\det(A))E,\\
r\overline{A}=-rA+\tr(A)rE.
\end{array}
$$
In the first formula, we used the Cayley--Hamilton theorem applied
to $\overline{A}$ and the identity $A+\overline{A}=\tr(A)E$.
$\Box$

We will now describe all sublattices of $\Mat_2(\Z)$ stable under
at least one of the operations $S_1,\dots, S_4$.
First consider sublattices of $\Mat_2(\Z)$ that are stable under
the multiplication, i.e. the operation $S_1$.

\begin{lemma}
Any non-null 2-dimensional sublattice of $\Mat_2(\Z)$ stable under $S_1$ contains
a nonzero scalar matrix.
\end{lemma}

\proof
Consider a non-null sublattice $\Lambda$ of $\Mat_2(\Z)$ stable under $S_1$.
Let $A$ be any matrix from $\Lambda$ such that $\det(A)\ne 0$.
Then the square $A^2$ also
belongs to $\Lambda$. By the Cayley--Hamilton theorem,
$A^2-\tr(A)A+\det(A)E=0$. Therefore, $\det(A)E\in\Lambda$.
$\Box$

\begin{prop}
\label{sublatS1}
Suppose that a non-null 2-dimensional sublattice of $\Mat_2(\Z)$ is
stable under $S_1$.
Then it is spanned by an integer matrix $A$ and an integer scalar matrix
$rE$, where $r\in\Z$ is an integer such that $\det(A)$ is divisible by $r$.
\end{prop}

\proof
Let $r$ be the smallest positive integer such that $rE\in\Lambda$.
Then, by Proposition \ref{basis}, there is a matrix $A\in\Lambda$
such that $A$ and $rE$ form a basis of $\Lambda$.
Since $\Lambda$ is stable under $S_1$, the matrix $A^2$
belongs to $\Lambda$. By the Cayley--Hamilton
theorem, $A^2-\tr(A)A=-\det(A)E\in\Lambda$, which means that
$r$ divides the determinant of $A$. $\Box$

The pairings $S_2$ and $S_3$ differ only by the order of arguments.
Hence it suffices to find all 2-dimensional sublattices of $\Mat_2(\Z)$
stable under $S_2$.

\begin{lemma}
If a non-null 2-dimensional sublattice of $\Mat_2(\Z)$ is stable under $S_2$, then
it contains a nonzero scalar matrix.
\end{lemma}

\proof
Consider a non-null sublattice $\Lambda$ of $\Mat_2(\Z)$ stable under $S_2$.
Let $A$ be any matrix from $\Lambda$ such that $\det(A)\ne 0$.
Then $\det(A)E=\overline{A}A$ also belongs to $\Lambda$.
$\Box$

\begin{prop}
\label{sublatS2}
Suppose that a non-null 2-dimensional sublattice of $\Mat_2(\Z)$ is
stable under $S_2$. Then it is spanned by an integer matrix $A$ and
an integer scalar matrix
$rE$, where $r\in\Z$ is an integer such that $\det(A)$ is divisible by $r$.
\end{prop}

\proof
Let $r$ be the smallest positive integer such that $rE\in\Lambda$.
Then, by Proposition \ref{basis}, there is a matrix $A\in\Lambda$
such that $A$ and $rE$ form a basis of $\Lambda$. Since
$\overline{A}A=\det(A)E\in\Lambda$, the determinant of $A$ is
divisible by $r$. $\Box$

The following proposition can be proved in exactly the same way
as Proposition \ref{sublatS2}.

\begin{prop}
\label{sublatS3}
Suppose that a non-null 2-dimensional sublattice of $\Mat_2(\Z)$ is
stable under $S_3$. Then it is spanned by an integer matrix $A$ and
an integer scalar matrix
$rE$, where $r\in\Z$ is an integer such that $\det(A)$ is divisible by $r$.
\end{prop}

Now consider sublattices of $\Mat_2(\Z)$ that are stable under $S_4$.

\begin{lemma}
\label{scalarS4}
Any 2-dimensional sublattice of $\Mat_2(\Z)$ stable under the pairing $S_4$ contains a
nonzero scalar matrix.
\end{lemma}

\proof
Suppose that a non-null sublattice $\Lambda$ of $SL(2,\Z)$ is stable
under the pairing $S_4$.

Take any matrix $A\in\Lambda$ such that $\tr(A)^2-\det(A)\ne 0$.
Such matrix exists because the quadratic map
$X\mapsto \tr(X)^2-\det(X)$ on $\Mat_2(\R)$ is positive definite on the 3-dimensional
subspace of $\Mat_2(\R)$ consisting of all symmetric matrices.
Hence this quadratic form can not vanish identically on the 2-dimensional subspace
spanned by $\Lambda$ over $\R$.

By the Cayley--Hamilton theorem, we have
$\overline{A}^2-\tr(A)\overline{A}+\det(A)E=0$.
We know that $\overline{A}^2\in\Lambda$. Therefore,
$$
\tr(A)A+\overline{A}^2=(\tr(A)^2-\det(A))E\in\Lambda.
$$
Here we used the identity $A+\overline{A}=\tr(A)E$.
$\Box$

\begin{prop}
\label{sublatS4}
Suppose that a 2-dimensional sublattice of $\Mat_2(\Z)$ is
stable under $S_4$. Then it is spanned by an integer matrix $A$ and
an integer scalar matrix
$rE$ such that $\tr(A)^2-\det(A)$ is divisible by $r$.
\end{prop}

\proof
Let $r$ be the smallest positive integer such that $rE\in\Lambda$.
Then, by Proposition \ref{basis}, there is a matrix $A\in\Lambda$
such that $A$ and $rE$ form a basis of $\Lambda$.
We know that $(\tr(A)^2-\det(A))E\in\Lambda$ (see the proof of
the preceding lemma). Therefore, $\tr(A)^2-\det(A)$ is divisible
by $r$.
$\Box$

The restrictions of the pairings $S_1$ through $S_4$ to non-null sublattices
stable under are given by formulas from Theorem \ref{inps} for some
values of the parameters.

\begin{thm}
\label{S1-32s1-3}
Let $A\in\Mat_2(\Z)$ be any integer $2\times 2$ matrix.
Take an integer $r$ that divides the determinant of $A$.
Consider the embedding $\phi:\R^2\to\Mat_2(\R)$ mapping a
vector with coordinates $(x_1,x_2)$ to the matrix $x_1A+x_2rE$.
Then, for $k=1$, $2$ or $3$, there is a unique integer normed pairing
$s_k:\R^2\times\R^2\to\R^2$ from Theorem \ref{inps} such that
$f_k(\x)=\det(\phi(\x))$ and
$$
\phi(s_k(\x,\y))=S_k(\phi(\x),\phi(\y)),\quad \x,\y\in\R^2.
$$
\end{thm}

\proof
Let us prove the theorem in the case $k=1$.
Other cases are analogous.
We have
$$
(x_1A+x_2rE)(y_1A+y_2rE)=x_1y_1A^2+(x_1y_2+x_2y_1)rA+x_2y_2r^2E=
$$
$$
\left(x_1y_1\tr(A)+(x_1y_2+x_2y_1)r\right)A+
\left(x_1y_1\frac{\det(A)}r+x_2y_2r\right)rE.
$$
Therefore, we should have
$$
s_1=\left(\begin{array}{cc}
\tr(A)& r\\ r& 0
\end{array}
\vline
\begin{array}{cc}
\frac{\det(A)}r& 0\\ 0& r
\end{array}\right)
$$
This matches the formula for $s_1$ from Theorem \ref{inps}, where we set
$$
p=0,\quad q=1,\quad m=-\frac{\det(A)}r,\quad k=\tr(A),\quad n=r
$$
The theorem is thus proved.
$\Box$

\begin{thm}
\label{S42s4}
Let $A$ be an integer $2\times 2$ matrix.
Take an integer $r$ that divides the number $\tr(A)^2-\det(A)$.
Consider the embedding $\phi:\R^2\to\Mat_2(\R)$ mapping a
vector with coordinates $(x_1,x_2)$ to the matrix $x_1A+x_2rE$.
Then there is a unique integer normed pairing
$s_4:\R^2\times\R^2\to\R^2$ from Theorem \ref{inps} such that
$f_4(\x)=\det(\phi(\x))$ and
$$
\phi(s_4(\x,\y))=S_4(\phi(\x),\phi(\y)),\quad \x,\y\in\R^2.
$$
\end{thm}

\proof
We have
$$
(x_1 \overline A+x_2 rE)(y_1 \overline A + y_2 rE)=z_1A+z_2 rE,
$$
where $z_1$ and $z_2$ are certain bilinear functions of $x=(x_1,x_2)$ and $y=(y_1,y_2)$.
Using the  relations
$\overline A^2=\tr(A)\overline A -\det(A)E$ and  $\overline A
=\tr(A)E-A$, we obtain:
$$
\begin{array}{l} z_1= -\tr(A)x_1y_1+ r (x_1 y_2+y_1x_2), \\
z_2=  \frac{\tr(A)^2-\det(A)}{r} x_1x_2+\tr(A)(x_1 y_2+y_1x_2) +r
x_2y_2.
\end{array}
$$
This matches the formula for $s_4$ from Theorem \ref{inps}, where
we set
$$
a=-\tr(A),\quad  b=0,\quad c=-r,\quad
d=-\frac{\tr(A)^2-\det(A)}r.
$$
The theorem is thus proved.
$\Box$

\section{Trigroup laws}

Let $f:\R^2\to\R$ be an integer quadratic form.
The product of two values of $f$ at points of the integer lattice is not necessarily
attained as a value at a point of the integer lattice.
In other words, not every integer quadratic form on $\R^2$ has
semigroup property.
But, as Vladimir Arnold discovered \cite{Ar}, the product of three values
$f(\x)f(\y)f(\z)$, where $\x,\y,\z\in\Z^2$, is always a value $f(\w)$
at some point $\w\in\Z^2$.
Moreover, there exists an integer 3-linear form
$[\cdot,\cdot,\cdot]:(\x,\y,\e)\mapsto [\x,\y,\e]$ such that
$$
f([\x,\y,\e])=f(\x)f(\y)f(\e).
$$
Such a 3-linear form is called a {\em trigroup law}.

There is an explicit formula for a 3-group law associated with a quadratic
form $f$.
Denote by $F$ the {\em polarization} of $f$, i.e.
a bilinear symmetric form such that $F(\x,\x)=f(\x)$ for all $\x\in\R^2$.
Then the formula
\begin{equation}
\label{3gr}
[\x,\y,\e]=-F(\x,\y)\e+F(\x,\e)\y+F(\y,\e)\x,\quad \x,\y,\e\in\R^2
\end{equation}
defines a trigroup law.
Formula (\ref{3gr}) was first found in \cite{Ai}.
In this section, we will give a new proof of the fact that formula (\ref{3gr})
provides a 3-group law with respect to the form $f$.
Other trigroup laws correspond to the permutations of
arguments and the change of the sign.
There are at most 6 trigroup laws that can be obtained from formula (\ref{3gr})
in this way.

Note that formula (\ref{3gr}) is asymmetric: the vector $\e$ plays a
distinguished role.

It is not obvious that the trilinear map (\ref{3gr})
is defined over integers, i.e. it takes any triple of integer vectors
to an integer vector.
In fact, for an integer quadratic form $f$,
the values of the polarization $F$ can be half-integers.
However, in formula (\ref{3gr}), all these half-integers
sum up miraculously to form an integer.

Let $\eps$ be $\pm 1$.
Denote by $\A$ the quotient algebra $\R[t]/(t^2+\eps)$.
Thus $\A$ is the field of complex numbers if $\eps=1$,
and $\A$ is the algebra of hyperbolic numbers if $\eps=-1$.
We introduced this unified notation to treat the cases of
positive definite and indefinite quadratic forms similarly.
The algebra $\A$ depends on $\eps$.
To emphasize this, we will sometimes write $\A(\eps)$ instead of $\A$.

Denote by $j$ the image of the polynomial $t$ under the
canonical projection of $\R[t]$ onto $\A$.
Thus we have $j^2=-\eps$.
{\em Conjugation} is the map of $\A$ to itself defined as follows:
$$
\overline{\xi_1+\xi_2 j}=\xi_1-\xi_2 j,\quad\xi_1,\xi_2\in\R.
$$
It is easy to verify that conjugation is an automorphism of $\A$.
Define the square of the norm by the formula: $|z|^2=z\overline{z}$.
Thus
$$
|\xi_1+\xi_2 j|^2=\xi_1^2+\eps\xi_2^2.
$$
We see that for any $z\in\A$, the number $|z|^2$ is a real number.
Note that $|z|$ is not always a real number, because for $\eps=-1$,
the number $|z|^2$ can be negative.

The algebra $\A$ is isomorphic to $\R^2$ as a real vector space.
Moreover, the basis $(1,j)$ defines an orientation in $\A$.
We will always assume $\A$ to be equipped with this standard orientation.
The coordinate space $\R^2$ also has the standard orientation given
by the standard basis $(e_1,e_2)$, where $e_1=(1,0)$ and $e_2=(0,1)$.

\begin{lemma}
For any nondegenerate quadratic form on $\R^2$ that is not negative definite,
there exists an orientation preserving real vector space isomorphism
$\phi:\R^2\to\A(\eps)$ such that $f(\x)=|\phi(\x)|^2$ for all $\x\in\R^2$.
Here $\eps=1$ for a positive definite form $f$, and $\eps=-1$ for an
indefinite form $f$.
\end{lemma}

\proof
This is just a reformulation of the well-known fact that any nondegenerate
quadratic form reduces to the sum of squares with coefficients $\pm 1$.
$\Box$

Two quadratic forms on $\R^2$ are said to be {\em properly equivalent}
if one can be obtained from the other by an integer invertible change
of coordinates.
In other words, a quadratic form $f$ is properly equivalent to a quadratic
form $f'$ if there is a matrix $A\in\SL_2(\Z)$ such that
$f'(\x)=f(A\x)$ for all $\x\in\R^2$.
A {\em class of forms} is a class of proper equivalence.

There is a natural correspondence between integer quadratic forms
on $\R^2$ and integer-normed lattices in $\A$.

Given a nondegenerate quadratic form $f$ that is not negative definite,
consider any orientation preserving real vector space
isomorphism $\phi:\R^2\to\A$ such that $f(\x)=|\phi(\x)|^2$ for all $\x\in\R^2$.

The image $L$ of $\Z^2$ is a certain two-dimensional lattice in $\A$.
It is {\em integer-normed}, i.e. for any element $z\in L$ the number $|z|^2$ is an integer.
Let us say that the lattice $L$ {\em corresponds} to the quadratic form $f$.

Define a {\em rotation} in $\A$ as an orientation preserving real vector
space automorphism $R:\A\to\A$ such that $|Rz|^2=|z|^2$ for all $z\in\A$.
All rotations of $\A$ form a group.

The following theorem is very classical:

\begin{thm}
\label{FormvsLat}
Choose $\eps=\pm 1$.
Classes of integer forms (positive definite if $\eps=1$, and indefinite if $\eps=-1$)
with respect to proper equivalence are in one-to-one correspondence with integer-normed lattices
in $\A$ considered up to a rotation.
\end{thm}

Felix Klein argued that Gauss knew this fact and used it to obtain his results.
But Gauss never mentioned this fact.

\proof
First note that a lattice corresponding to a given quadratic form
is uniquely defined up to a rotation.

Let integer quadratic forms $f,f':\R^2\to\R$ be properly equivalent.
Then there is a map $A:\R^2\to\R^2$ such that
$A$ takes the lattice $\Z^2$ to itself, preserves the
orientation and satisfies the property $f'=f\circ A$.

Let $\phi:\R^2\to\A$ be an orientation preserving
isomorphism of real vector spaces such that $f=|\phi|^2$.
Then $f'=|\phi\circ A|$, and $\phi\circ A$ takes
the lattice $\Z^2$ to the same lattice as $\phi$.
Thus an integer quadratic form determines an integer-normed
lattice uniquely up to a rotation.

Conversely, given an integer-normed lattice $L$ in $\A$,
we can recover a corresponding quadratic form up to
proper equivalence.
Namely, let us choose a basis $e_1,e_2$ in $L$.
Suppose that this basis is {\em properly oriented}, i.e.
it provides the same orientation of $\A$ as the standard basis $(1,j)$.
The choice of a properly oriented basis gives rise to an orientation preserving
isomorphism $\phi:\R^2\to\A$ defined by the formula $\phi(\x)=x_1e_1+x_2e_2$ for
any point $\x\in\R^2$ with coordinates $(x_1,x_2)$.
Define the quadratic form $f$ to be $f=|\phi|^2$.
This quadratic form depends on the particular choice
of a basis in $L$, but different bases give rise to properly equivalent forms.
$\Box$

We will now proceed to the proof that formula (\ref{3gr}) provides a trigroup law.

\begin{thm}
\label{3grinter}
Consider 3 vectors $\x,\y,\e\in\R^2$. If $\phi(\x)=u$, $\phi(\y)=v$ and
$\phi(\e)=w$, then
$$
\phi[\x,\y,\e]=uv\overline{w},
$$
where $[\x,\y,\e]$ is the value of the trigroup law given by formula
(\ref{3gr}).
\end{thm}

\proof
Since the quadratic form $f$ corresponds to the square of the absolute value
under the isomorphism $\phi$, we have
$$
\phi\circ F(\x,\y)=\frac12(u\overline{v}+\overline{u}v),\quad
\phi\circ F(\x,\e)=\frac12(u\overline{w}+\overline{u}w),\quad
\phi\circ F(\y,\e)=\frac12(v\overline{w}+\overline{v}w).
$$
We can now plug in these expressions to formula (\ref{3gr}):
$$
\phi[\x,\y,\e]=-\frac12(u\overline{v}+\overline{u}v)w+\frac12(u\overline{w}+\overline{u}w)v
+\frac12(v\overline{w}+\overline{v}w)u=uv\overline{w}.
$$
$\Box$

\begin{cor}
For any triple of vectors $\x,\y,\e\in\R^2$ we have
$$
f([\x,\y,\e])=f(\x)f(\y)f(\e).
$$
\end{cor}

\proof
Indeed, the quadratic form $f$ corresponds to the square of the absolute
value under the isomorphism $\phi$.
Set $u=\phi(\x)$, $v=\phi(\y)$, $w=\phi(\e)$. We have
$$
f([\x,\y,\e])=|uv\overline{w}|^2=|u|^2|v|^2|w|^2=f(\x)f(\y)f(\e).
$$
$\Box$

Let us now prove that the trigroup law defined by formula (\ref{3gr})
takes triples of integer vectors to integer vectors.
This statement can be reformulated as follows:

\begin{prop}
Suppose that a 2-dimensional lattice $L\subset\A$ is such that
for any $z\in L$ the number $|z|^2$ is integer. Then $L$ is
stable under the 3-linear operation $(u,v,w)\mapsto uv\overline{w}$.
In other words, for any triple $u,v,w\in L$, we have
$uv\overline{w}\in L$.
\end{prop}

\proof
Choose any basis $(z_1,z_2)$ of $L$. It suffices to show that $z_1^2\overline{z_2}$
is an integer linear combination of $z_1$ and $z_2$.
Since $|z_1+z_2|^2$ is an integer, we have that $a=z_1\overline{z_2}+\overline{z_1}z_2$
is also an integer. Then
$$
z_1^2\overline{z_2}=z_1(a-\overline{z_1}z_2)=az_1-|z_1|^2z_2.
$$
$\Box$

We will need the following proposition in the sequel:

\begin{prop}
\label{conf}
Suppose that a map $\Phi:\A\to\A$ is conformal, i.e. there is
a constant $C\ne 0$ such that $|\Phi(w)|^2=C|w|^2$ for any $w\in\A$.
Then there exists an element $z\in\A$ such that
$\Phi$ is either the multiplication by $z$, i.e. $\Phi(w)=zw$ for all $w\in\A$,
or the composition of multiplication by $z$ and the conjugation,
i.e. $\Phi(w)=z\overline{w}$ for all $w\in\A$.
\end{prop}

\proof
Set $z=\Phi(1)$. Since $C\ne 0$, the element $z$ is invertible.
Consider the map $w\mapsto z^{-1}\Phi(w)$. It fixes 1 and
preserves the norms.
But there are only two such maps, namely,
the identity and the automorphism of conjugation.
$\Box$

\begin{remark}
For a negative form $f$, the form $-f$ is positive definite,
and any trigroup law for $-f$ is also a trigroup law for $f$.
\end{remark}

The following remarks are due to V. Arnold.
Suppose that an integer quad\-ra\-tic form $f$ attains value 1 at some
point $\e_0\in\Z^2$.
Then any trigroup law gives rise to a normed
pairing with respect to $f$.
Namely, if $(\x,\y,\e)\mapsto[\x,\y,\e]$ is a trigroup law, then
$(\x,\y)\mapsto [\x,\y,\e_0]$ is a normed pairing.

A more general situation is as follows:

\begin{thm}
\label{3gr2inp}
Suppose that a quadratic form $f$ equals to $rg$, where $r$ is an integer,
and $g$ is an integer quadratic form attaining the value $r$.
Denote by $\e_0$ any integer vector such that $g(\e_0)=r$.
Then, for any trigroup law $(\x,\y,\e)\mapsto[\x,\y,\e]$ for the form $f$,
the pairing $(\x,\y)\mapsto [\x,\y,\e_0]/r$ is a normed pairing with respect to $f$.
\end{thm}

\proof
Set $\z=[\x,\y,\e_0]/r$.
We have
$$
f(\z)=f(\x)f(\y)f(\e_0)/r^2=f(\x)f(\y).
$$
$\Box$

Theorem \ref{3gr2inp} gives an alternative construction of integer normed pairings.
Let us explain the relation between this construction and Theorem \ref{inps}.

\begin{thm}
\label{3gr2s1-3}
Let $g=(m,k,n)$ be an arbitrary integer quadratic form.
Denote by $G$ its polarization.
Suppose that $[\cdot,\cdot,\cdot]$ is the trigroup law for $f$ given
by formula (\ref{3gr}).
Take any vector $\e\in\Z^2$ with coordinates $(p,q)$.
Pairings $s_1$ through $s_3$ from Theorem \ref{inps} coincide
with the following pairings:
$$
\begin{array}{l}
s_1(\x,\y) = -G(\x,\y)\e + G(\x,\e)\y + G(\y,\e)\x=[\x,\y,\e]/r\\
s_2(\x,\y) = G(\x,\y)\e - G(\y,\e)\x + G(\x,\e)\y=[\y,\e,\x]/r\\
s_3(\x,\y) = G(\x,\y)\e - G(\x,\e)\y + G(\y,\e)\x=[\x,\e,\y]/r.
\end{array}
$$
The parameters $m,k,n,p,q$ from Theorem \ref{inps} match
the numbers $m,k,n,p,q$ introduced above.
The pairings $s_1$ through $s_3$ are integer-normed pairings
with respect to the quadratic form $f=rg$, where $r=g(\e)$.
\end{thm}

This theorem can be proved by a direct computation.

\section{Normed pairings and lattices I}

In this section and in the next one, we describe integer normed pairings from
Theorem \ref{inps} in terms of lattices in $\A$.

Fix a normed pairing $s:\R^2\times\R^2\to\R^2$ with respect to a nondegenerate
quadratic form $f:\R^2\to\R$.
Recall that $s$ is a {\em normed pairing} if $f(s(\x,\y))=f(\x)f(\y)$ for all
$\x,\y\in\R^2$.
For the time being, we do not require that $s$ be defined over integers.
Note that the form $f$ can not be negative definite.
Otherwise, $f(s(\x,\y))=f(\x)f(\y)$ is positive as the product of two negative
numbers, which is a contradiction.

Set $\eps=1$, if the form $f$ is positive definite, and $\eps=-1$, if
$f$ is indefinite.
Fix an orientation preserving real vector space isomorphism
$\phi:\R^2\to\A$ such that $f(\x)=|\phi(\x)|^2$ for all $\x\in\R^2$.
Let us introduce the $\R$-bilinear pairing $\sigma:\A\times\A\to\A$ defined
by the formula $\sigma(z,w)=\phi(s(\phi^{-1}(z),\phi^{-1}(w))$ for all $z,w\in\A$.
Then $\sigma$ is a normed pairing with respect to the square of the absolute value,
i.e. $|\sigma(z,w)|^2=|z|^2|w|^2$ for all $z,w\in\A$.

Recall the proof of the following classical theorem:

\begin{thm}
\label{Hu}
There exists a number $\alpha\in\A$ such that $|\alpha|^2=1$ and
$\sigma=\alpha\sigma_k$ for some $k=1,\dots, 4$, where $\sigma_k$ are the
pairings defined in Theorem \ref{sigma1-4}.
\end{thm}

\proof
Fix a point $z\in\A$ and consider the map $\sigma_z:w\mapsto \sigma(z,w)$.
This map is conformal (it multiplies all distances by $|z|$).
Therefore, it has either the form $\sigma_z:w\mapsto A(z)w$ or the form
$\sigma_z:w\mapsto A(z)\overline{w}$
where $A(z)$ is an element of $\A$ such that $|A(z)|^2=|z|^2$.
This follows from Proposition \ref{conf}.

Note that $A(z)$ depends linearly on $z$.
Hence, we have a linear map $A$ preserving the distances.
It follows that $A:z\mapsto \alpha z$ or $A:z\mapsto \alpha\overline{z}$,
where $\alpha\in\A$ is such that $|\alpha|^2=1$.
The theorem now follows.
$\Box$

Actually, we can even get rid of the number $\alpha$.
To this end, we need the following lemma:

\begin{lemma}
Any element $\alpha\in\A$ such that $|\alpha|^2=1$ has a cubic root.
In other words, there exists an element $\beta\in\A$ such that $\beta^3=\alpha$.
\end{lemma}

\proof
If $\A=\C$, this is obvious.
Suppose that $\A=\H$, the algebra of hyperbolic numbers.

Note that for every real number $t$, we have
$e^{jt}=\cosh(t)+j\sinh(t)$, and that any number
$\alpha$ such that $|\alpha|^2=1$ has the form $\pm e^{jt}$.
It suffices to take $\beta=\pm e^{jt/3}$.
$\Box$

\begin{prop}
\label{alpha1}
By changing the isomorphism $\phi$, we can arrange that $\alpha=1$.
\end{prop}

\proof
Consider an isomorphism $\phi':\R^2\to\A$ that maps any vector $x\in\R^2$
to $\alpha'\phi(x)$.
Here $\alpha'$ is some element of $\A$ such that $|\alpha'|^2=1$.
We will fix this element later.
Denote by $\sigma'$ the normed pairing that corresponds
to $s$ under the isomorphism $\phi'$. We have
$$
\sigma'(z,w)=\alpha'^{-1}\sigma(\alpha'z,\alpha'w).
$$

The proper choice of the number $\alpha'$ depends on what operation
$\sigma$ we have.
If $\sigma=\alpha\sigma_1$, then choose $\alpha'=\alpha^{-1}$ so that $\sigma'=\sigma_1$.
If $\sigma=\alpha\sigma_2$ or $\alpha\sigma_3$, then choose $\alpha'=\alpha$ so that $\sigma'=\sigma_2$ or $\sigma_3$.
Finally, if $\sigma=\alpha\sigma_4$, then choose $\alpha'=\sqrt[3]{\alpha}$ so that $\sigma'=\sigma_4$.
By the previous lemma, the cubic root exists.
$\Box$

\label{p:sigma1-4}
\subsection*{Proof of Theorem \ref{sigma1-4}.}
Let $f$ be a nondegenerate quadratic form.
Suppose that there exists an integer normed pairing $s$ with respect to $f$.
As we saw, $f$ must be either positive definite or indefinite in this case.
Set $\eps=1$, if $f$ is positive definite, and $\eps=-1$, if $f$ is indefinite.
Consider an orientation preserving vector space isomorphism $\phi:\R^2\to\A(\eps)$.
Let $L$ be the image of $\Z^2$ under the isomorphism $\phi$.

By Theorem \ref{Hu}, the lattice $L$ is stable under a pairing $\alpha\sigma_k$,
where $k=1,2,3$ and $\alpha$ is an element of $\A=\A(\eps)$ such that $|\alpha|^2=1$.
By Proposition \ref{alpha1}, we can now assume without loss of generality that $\alpha=1$.
Otherwise we just need to change the isomorphism $\phi$.
Thus we have proved Theorem \ref{sigma1-4}.
$\Box$

\begin{thm}
\label{s1-32S1-3}
Suppose that a 2-dimensional lattice $L\subset\A$ is stable under $\sigma_k$
for some $k=1,2,3$.
Then there exists an embedding $\iota: L\to \Mat_2(\Z)$ such that
$\det(\iota(z))=|z|^2$ for all $z\in L$, and the lattice $\iota(L)$ is stable
under $S_k$.
Moreover, $\sigma_k$ corresponds to $S_k$ under this embedding:
$$
\iota\circ\sigma_k(z,w)=S_k(\iota(z),\iota(w)),\quad z,w\in\A
$$
\end{thm}

\proof
{\em Case $k=1$.}
For any $z\in L$, consider the operator $\iota(z)$ of multiplication by $z$.
This operator takes the lattice $L$ to itself.
Choose a basis in $L$.
The matrix of the operator $\iota(z)$ in this basis is an integer $2\times 2$ matrix.
Slightly abusing notation, denote this matrix by $\iota(z)$, the same as the corresponding
operator.

The determinant of $\iota(z)$ is equal to $|z|^2$.
Indeed, $\iota(z)$ is a conformal operator: $|\iota(z)w|^2=|z|^2|w|^2$ for all $w\in\A$.
It is obvious that $\iota(zw)=\iota(z)\iota(w)$ for all $z,w\in L$.

{\em Case $k=2$.}
For any number $z\in L$, define the operator $\iota(z)$ by the formula
$\iota(z):w\mapsto \overline{z}w$.
This operator takes the lattice $L$ to itself.
Analogously to the case $k=1$, the determinant of $\iota(z)$ is $|z|^2$.
The adjoint operator to $\iota(z)$ is clearly $\iota(\overline{z})$.
Therefore,
$\iota(\overline{z}w)=\overline{\iota(z)}\iota(w)$ for all $z,w\in L$.

{\em Case $k=3$} is completely analogous to the case $k=2$.
$\Box$

Let us now give a geometric description of lattices stable under
one of the operations $\sigma_k$, $k=1,2,3$.
An element $z\in\A-\R$ is called a {\em quadratic integer} if it satisfies
a quadratic equation $z^2+bz+c=0$ with integer coefficients $b,c\in\Z$.
Note that our terminology is different from the usual one,
according to which a real number can be a quadratic integer.
Even if the equation $z^2+bz+c=0$ has real roots, we are
interested in imaginary (non-real) roots only.

\begin{prop}
Any quadratic equation over $\R$, whose discriminant is nonzero,
has exactly 2 imaginary roots in $\A$.
They are conjugate to each other.
\end{prop}

\proof
Indeed, it is easy to see that imaginary roots of a quadratic equation
$az^2+bz+c=0$ are given by the formula
$$
z_{1,2}=\frac{-b\pm j\sqrt{|b^2-4ac|}}{2a}.
$$
$\Box$

The following corollary is an analog of the Vieta theorem:

\begin{cor}
\label{Vieta}
An element $\zeta\in\A$ is a quadratic integer if and only if
the numbers $\zeta+\overline{\zeta}$ and $|\zeta|^2$ are integers.
\end{cor}

\begin{thm}
\label{zL}
Suppose that a lattice $L\subset\A$ is stable under multiplication by $z\in\A$,
i.e. $zw\in L$ for all $w\in L$.
Then $z$ is a quadratic integer.
\end{thm}

\proof
Consider the operator $Z$ of multiplication by $z$.
Choose a basis for the lattice $L$ and denote by $\hat Z$ the matrix of
$Z$ in this basis.
Since $Z$ preserves the lattice $L$, the matrix $\hat Z$ has
integer entries.
Therefore, the trace and the determinant of $Z$ are integers.
Set $b=-\tr(Z)$, $c=\det(Z)$.
By the Cayley--Hamilton theorem, $Z$ satisfies the equation $Z^2+bZ+c=0$.
The theorem now follows.
$\Box$

\begin{lemma}
\label{lat-s1-sq}
Suppose that a lattice $L\subset\A$ is stable under multiplication
(i.e. $L$ is stable under operation $\sigma_1$).
Then for any $z\in L$, the number $|z|^2$ also belongs to $L$.
\end{lemma}

\proof
By Theorem \ref{zL}, any element $z\in L$ satisfies a
quadratic equation $z^2+bz+c=0$ with integer coefficients $b$ and $c$.
The vector $bz$ is in $L$.
Since $L$ is stable under multiplication, $z^2$ also belongs to $L$.
Therefore, $c\in L$.
But $c=|z|^2$ by the Vieta theorem.
$\Box$

\begin{thm}
\label{lat-sigma1}
Suppose that a lattice $L\subset\A$ is stable under $\sigma_1$.
Then $L$ is generated by an integer $r\in\Z$ and a quadratic integer
$\zeta\in\A$ such that $|\zeta|^2$ is divisible by $r$.
In particular, the lattice $L$ is integer-normed.
\end{thm}

\proof
By Lemma \ref{lat-s1-sq}, the lattice $L$ contains some
nonzero integers.

Denote by $r$ the smallest positive integer contained in $L$.
By Proposition \ref{basis}, there exists a basis in $L$ of the form $(r,\zeta)$.
The lattice $L$ is stable under $\sigma_1$.
In particular, it is stable under multiplication by $\zeta$.
From Theorem \ref{zL}, it follows that $\zeta$ is a quadratic integer.
We know that $|\zeta|^2\in L$.
Therefore, $|\zeta|^2$ is divisible by $r$.

For any pair of integers $m_1,m_2\in\Z$, we have
$$
|m_1r+m_2\zeta|^2=m_1^2r^2+m_1m_2r(\zeta+\overline{\zeta})+m_2^2|\zeta|^2\in\Z
$$
because $\zeta+\overline{\zeta}$ and $|\zeta|^2$ are integers by
the Vieta theorem.
$\Box$

\begin{thm}
\label{lat-sigma23}
Suppose that a lattice $L\subset\A$ is stable under $\sigma_2$ or $\sigma_3$.
Then $L$ is generated by an integer $r\in\Z$ and
a quadratic integer $\zeta\in\A$ such that
$|\zeta|^2$ is divisible by $r$.
In particular, the lattice $L$ is integer-normed.
\end{thm}

\proof
For any $z\in L$, we have $|z|^2=z\overline{z}\in L$.
Thus the lattice $L$ contains some nonzero integers.

Denote by $r$ the smallest positive integer contained in $L$.
By Proposition \ref{basis}, there exists a basis in $L$ of the form $(r,\zeta)$.
The lattice $L$ is stable under $\sigma_2$ or $\sigma_3$.
In particular, it is stable under multiplication by $\overline{\zeta}$.
From Theorem \ref{zL}, it follows that $\overline{\zeta}$ is a quadratic integer,
hence $\zeta$ is a quadratic integer as well.
We know that $|\zeta|^2\in L$.
Therefore, $|\zeta|^2$ is divisible by $r$.
$\Box$

\subsection*{Discriminants of quadratic forms and of lattices.}
Recall the definition of the discriminant and some of its basic properties.
The {\em discriminant} of a quadratic form $(m,k,n)$ is defined as
\begin{equation}
\label{discr-form}
\Delta=k^2-4mn.
\end{equation}
For positive definite forms, the discriminant is negative.
For indefinite forms, the discriminant is positive.

It is easy to see that the discriminant of any integer quadratic form on
$\R^2$ is congruent either to 0 or to 1 modulo 4.
Moreover, the discriminant of $(m,k,n)$ is divisible by 4 exactly when $k$
is even.

Now consider a two-dimensional lattice $L$ in $\A$.
Take any basis $e_1$, $e_2$ of $L$.
Then the {\em discriminant} of $L$ is defined as
\begin{equation}
\label{discr-lat}
\Delta=-4\det\left(\begin{array}{cc}
\<e_1,e_1\>& \<e_1,e_2\>\\
\<e_1,e_2\>& \<e_2,e_2\>
\end{array}\right)
\end{equation}
Here the bilinear form $(z,w)\mapsto \<z,w\>$ is the polarization of
the quadratic form $z\mapsto |z|^2$.

The factor of $-4$ is introduced here to match the notation for the
discriminant of a quadratic form.
It is clear that the discriminant does not depend on the particular
choice of a basis in $L$.
It is easy to see that
\begin{equation}
\label{discr}
\Delta=(e_1\overline{e_2}-e_2\overline{e_1})^2
\end{equation}

If $\A=\C$, then $\Delta$ is always a negative real number.
If $\A=\H$, then $\Delta$ is always positive.
Inner products of vectors from an integer-normed lattice are integers or half-integers.
Therefore, the discriminant of an integer-normed lattice is always a negative integer.
The discriminant $\Delta$ of an integer-normed lattice coincides with the
discriminant of any corresponding quadratic form.
In particular, the discriminants of properly equivalent forms coincide.

Let $L$ be a lattice in $\C$.
Choose a basis $(e_1,e_2)$ in $L$.
Denote by $\Pi$ the {\em fundamental parallelogram of $L$}, i.e.
the parallelogram spanned by the vectors $e_1$ and $e_2$.
The vertices of the parallelogram $\Pi$ are 0, $e_1$, $e_2$ and $e_1+e_2$.
The geometric meaning of the discriminant $\Delta$ of $L$ is that $\Delta=-4Area^2(\Pi)$.

\subsection*{The ring $R_\Delta$.}
Consider an integer $\Delta$ that is either divisible by 4 or congruent to 1 modulo 4.
Recall that for $\zeta\in\A$, the ring $\Z[\zeta]$ is defined as a subring of $\A$
generated by $\Z$ and $\zeta$.
Let us define the ring $R_\Delta$ as follows.
If $\Delta$ is divisible by 4, then $R_\Delta$ is $\Z[j\sqrt{|\Delta|}/2]$.
If $\Delta$ is congruent to 1 modulo 4, then $R_\Delta$ is $\Z[(1+j\sqrt{|\Delta|})/2]$.

Thus $R_\Delta$ is always a subring of $\A$.
Even if $\Delta>0$, the ring $R_\Delta$ does not lie in $\R$.

\begin{thm}
\label{rings}
Suppose that $\zeta\in\A$ is a quadratic integer.
Let $\zeta^2+b\zeta+c=0$ be a quadratic equation on $\zeta$.
Denote by $\Delta$ the discriminant of this equation: $\Delta=b^2-4c$.
Then $\Z[\zeta]$ is the same as $R_\Delta$.
\end{thm}

\proof
This follows from the explicit formula for the roots of a quadratic equation:
$$
\zeta=\frac{-b\pm j\sqrt{|\Delta|}}2.
$$
If $b$ is even, then $\zeta$ differs from $j\sqrt{|\Delta|}/2$ by an integer.
Therefore, $\Z[\zeta]=\Z[j\sqrt{|\Delta|}/2]$.
If $b$ is odd, then $\zeta$ differs from $(1+j\sqrt{|\Delta|})/2$ by an integer.
Therefore, $\Z[\zeta]=\Z[(1+j\sqrt{|\Delta|})/2]$.
We see that in both cases $\Z[\zeta]=R_\Delta$.
$\Box$

\begin{prop}
\label{R_Delta-iqi}
All elements of the ring $R_\Delta$ are quadratic integers.
\end{prop}

\proof
If $\Delta$ is divisible by 4, then any element $\zeta$ of $R_\Delta$ has
the form $m_1+m_2j\sqrt{|\Delta|}/2$, where $m_1$ and $m_2$ are integers.
Then $|\zeta|^2=m_1^2-m_2^2\Delta/4$ and $\zeta+\overline{\zeta}=2m_1$ are integers.
By Corollary \ref{Vieta}, it follows that $\zeta$ is a quadratic integer.

If $\Delta$ is congruent to 1 modulo 4, then any element $z$ of $R_\Delta$ has
the form $m_1+m_2(1+j\sqrt{|\Delta|})/2$, where $m_1$ and $m_2$ are integers.
Then
$|z|^2=m_1^2+m_1m_2-m_2^2(\Delta-1)/4$ and
$\zeta+\overline{\zeta}=2m_1+m_2$ are integers.
By Corollary \ref{Vieta}, it follows that $\zeta$ is a quadratic integer.
$\Box$

\begin{prop}
Suppose that a lattice $L\subset\A$ is generated by an integer $r$ and a
quadratic integer $\zeta\in\A$ such that $|\zeta|^2$ is divisible by $r$.
Then the lattice $L$ is stable under multiplication by $\zeta$, i.e.
$\zeta L\subseteq L$.
\end{prop}

\proof
It suffices to show that $\zeta^2\in L$.
Since $\zeta$ is a quadratic integer, we have $\zeta^2+b\zeta+rc=0$, where
$b$ is an integer and $c=|\zeta|^2/r$ is also an integer.
Thus $\zeta^2=-b\zeta-rc\in L$.
$\Box$

\begin{thm}
\label{ideal}
Suppose that a lattice $L\subset\A$ is generated by an integer $r$ and a
quadratic integer $\zeta\in\A$ such that $|\zeta|^2$ is divisible by $r$.
Denote by $\Delta$ the discriminant of $L$ and set $\Delta^*=\Delta/r^2$
(note that $\Delta^*$ must be an integer congruent to 0 or 1 modulo 4).
Then $L$ is an ideal in the ring $R_{\Delta^*}$.
\end{thm}

\proof
By formula (\ref{discr}), we have
$$
\Delta=(r\overline{\zeta}-r\zeta)^2=r^2(\zeta-\overline{\zeta})^2
$$
We can assume that $\zeta$ is in the upper half-plane, i.e. the basis
$(r,\zeta)$ is properly oriented.
Otherwise, just replace $\zeta$ with $\overline{\zeta}$.

We have:
$$
\zeta-\overline{\zeta}=\frac{j\sqrt{|\Delta|}}r,\quad \zeta+\overline{\zeta}=-b\in\Z.
$$
The second equality follows from the fact that $\zeta$ is a quadratic integer.
Solving this system of linear equations for $\zeta$, we obtain:
$$
\zeta=\frac{-br+j\sqrt{|\Delta^*|}}2
$$
It follows that $\zeta\in R_{\Delta^*}$.
Moreover, we have $R_{\Delta^*}=\Z[\zeta]$.
Hence, $L\subseteq R_{\Delta^*}$.

Let us now prove that $L$ is an ideal in $R_{\Delta^*}$.
Since $R_{\Delta^*}=\Z[\zeta]$, it is enough to verify that $L$ is stable under
multiplication by $\zeta$.
But this is exactly the statement of the previous proposition.
$\Box$

\label{p:sigma1-3}
\subsection*{Proof of Theorem \ref{sigma1-3}.}
Let $L$ be a lattice stable under one of operations $\sigma_k$, $k=1,2,3$.
From Theorems \ref{lat-sigma1} and \ref{lat-sigma23}, it follows that
$L$ is generated by an integer $r\in\Z$ and a quadratic integer $\zeta\in\A$
such that $|\zeta|^2$ is divisible by $r$.
In particular, the lattice $L$ is integer-normed.
Let $\Delta$ denote the discriminant of the lattice $L$.
It is divisible by $r^2$.
Set $\Delta^*=\Delta/r^2$.
Note that $\Delta^*$ is congruent to 0 or 1 modulo 4, since both
$\Delta$ and $r^2$ satisfy this property.
By Theorem \ref{ideal}, the lattice $L$ is an ideal in the ring $R_{\Delta^*}$.
By definition, the ring $R_{\Delta^*}$ is stable under conjugation.
Therefore, $\overline{L}\subset R_\Delta$, and we have
$LL\subseteq R_{\Delta^*} L=L$ and $\overline{L}L\subseteq R_{\Delta^*} L=L$.
This means that $L$ is stable under all operations $\sigma_k$, $k=1,2,3$.
Theorem \ref{sigma1-3} is thus proved.
$\Box$

\section{Normed pairings and lattices II}

In this section, we will describe lattices in $\A$ stable under $\sigma_4$.
In particular, we will prove the following theorem:

\begin{thm}
\label{lat-sigma4}
Suppose that a lattice $L\subset\A$ is stable under the
pairing $\sigma_4$.
Consider an orientation preserving real vector space isomorphism
$\phi:\R^2\to\A$ such that $L=\phi(\Z^2)$.
Define a quadratic form $f$ by the formula $f(\x)=|\phi(\x)|^2$
for all $\x\in\R^2$.
Set
$$
s(\x,\y)=\phi^{-1}\sigma_4(\phi(\x),\phi(\y)),\quad \x,\y\in\R^2
$$
Then $s$ is of type $(-,-)$.
\end{thm}

We will need two lemmas:

\begin{lemma}
\label{s-sc}
The pairing $s$ satisfies the following relation:
$$
s(\x,s(\x,\y))=f(\x)\y
$$
for all $\x,\y\in\R^2$.
\end{lemma}

\proof
Indeed, set $z=\phi(\x)$ and $w=\phi(\y)$.
Then
$$
\sigma(z,\sigma(z,w))=\overline{z\overline{zw}}=\overline{z}zw=|z|^2w.
$$
Applying $\phi^{-1}$ to the both parts of this equation,
we obtain the desired result.
$\Box$

\begin{lemma}
\label{s-trl}
Consider a linear operator $A:\R^2\to\R^2$ that is not a scalar operator.
The operator $A^2$ is a scalar operator if and only if the operator $A$ is traceless.
\end{lemma}

\proof
Indeed, by the Cayley--Hamilton theorem, we have $A^2-\tr(A)A+\det(A)=0$.
The lemma follows immediately from this relation.
$\Box$

\subsection*{Proof of Theorem \ref{lat-sigma4}.}
Let us show that the pairing $s$ is commutative and traceless.
Commutativity of $s$ follows immediately from commutativity of $\sigma_4$.
For $\x\in\R^2$, consider the linear operator $M_\x:\y\mapsto s(\x,\y)$.
We need to prove that $\tr(M_\x)=0$.
By Lemma \ref{s-sc}, the operator $M_\x^2$ is the scalar operator $f(\x)E$.
By Lemma \ref{s-trl}, the operator $M_\x$ must be traceless
(note that $M_\x$ can not be scalar because it reverses orientation).

Since $s$ is commutative and traceless, it is given by formula (\ref{cotrl}).
It follows that $s$ coincides with the pairing $s_4$ from Theorem \ref{inps}.
Moreover, the parameters $a$, $b$, $c$, $d$ from formula (\ref{cotrl}) match
the parameters $a$, $b$, $c$, $d$ from Theorem \ref{inps}.
Since $s$ is defined over integers, i.e. $s(\Z\times\Z)\subseteq\Z$,
the parameters $a$, $b$, $c$, $d$ must be integers.
$\Box$

\label{p:inps}
\subsection*{Proof of Theorem \ref{inps}.}
Suppose that $f$ is a nondegenerate quadratic form on $\R^2$, and that $s$ is
an integer normed pairing with respect to $f$.
In Section 4, we saw that the form $f$ corresponds to a lattice $L\subseteq\A$ that
is stable under one of the operations $\sigma_k$, $k=1,\dots, 4$.
Namely, for some orientation preserving vector space isomorphism
$\phi:\R^2\to\A$, we have $f=|\phi|^2$ and $L=\phi(\Z^2)$.

{\em Case $k=1,2,3$.}
First assume that $k=1,2,3$.
By Theorems \ref{lat-sigma1} and \ref{lat-sigma23}, the lattice $L$ is generated
by an integer $r$ and a quadratic integer $\zeta\in\A$ such that $|\zeta|^2$
is divisible by $r$.
Any element of $L$ has the form $m_1r+m_2\zeta$, where $m_1$ and $m_2$ are
integers.
We have
$$
|m_1r+m_2\zeta|^2=m_1^2r^2+m_1m_2r(\zeta+\overline{\zeta})+m_2^2|\zeta|^2.
$$
Since $|\zeta|^2$ is divisible by $r$, the numbers $|z|^2$ are divisible by
$r$ for all $z\in L$.
Therefore, the quadratic form $f$ has the form $rg$, where $g$ is another
integer quadratic form.

Note that $g$ attains the value $r$ at some integer vector.
Indeed, if $\e=\phi^{-1}(r)\in\Z^2$, then $f(\e)=r^2$, hence $g(\e)=r$.

The trilinear map $(u,v,w)\mapsto uv\overline{w}$ corresponds to a trigroup law
$(\x,\y,\e)\mapsto [\x,\y,\e]$ for the form $f$.
Namely, if $u=\phi(\x)$, $v=\phi(\y)$ and $w=\phi(\e)$, then
$uv\overline{w}=\phi[\x,\y,\e]$.
Moreover, this trigroup law is given by formula (\ref{3gr}).
This is proved in Proposition \ref{3grinter}.
Since all coefficients of $f$ are divisible by $r$, vectors $[\x,\y,\e]/r$
are always integer.
This means that the lattice $L$ is stable under the trilinear operation
$(u,v,w)\mapsto uv\overline{w}/r$.

For any fixed $w\in L$, consider the following bilinear operations:
$$
\begin{array}{c}
\sigma'_1:(u,v)\mapsto uv\overline{w}/r,\\
\sigma'_2:(u,v)\mapsto vw\overline{u}/r,\\
\sigma'_3:(u,v)\mapsto uw\overline{v}/r,
\end{array}
$$
The lattice $L$ is stable under all of these operations.
If we set $w=r$, then $\sigma'_k=\sigma_k$ for all $k=1,2,3$.

On the other hand, operations $\sigma'_k$ correspond to integer
normed pairings $s_k$ under the isomorphism $\phi$.
Namely,
$$
\begin{array}{c}
s_1(\x,\y)=[\x,\y,\e]/r=\phi^{-1}(uv\overline{w}/r),\\
s_2(\x,\y)=[\y,\e,\x]/r=\phi^{-1}(vw\overline{u}/r),\\
s_3(\x,\y)=[\x,\e,\y]/r=\phi^{-1}(uw\overline{v}/r),
\end{array}
$$
where $u=\phi(\x)$, $v=\phi(\y)$ and $w=\phi(\e)$.
This follows from Theorem \ref{3gr2s1-3}.
Thus, in particular, $s$ is of type $(+,+)$, $(-,+)$ or $(-,-)$.

{\em Case $k=4$.}
Now assume that $k=4$.
Then the lattice $L$ is stable under $\sigma_4$.
By Theorem \ref{lat-sigma4}, it follows that the pairing $s$ coincides
with $s_4$.
$\Box$

\section{Examples}

We saw that integer normed pairings of types $(+,+)$, $(-,+)$ and $(+,-)$ always
come together: if a quadratic form admits an integer normed
pairing of one type, then the same form admits integer normed pairings of two other
types.

However, there exist quadratic forms that admit integer normed pairings of types
$(+,+)$, $(-,+)$ and $(+,-)$ but do not admit integer normed pairings of type
$(-,-)$.
Below, we will give an example of such a quadratic form.

Let $f$ be a quadratic form that is either positive definite, or indefinite.
If $f$ is positive definite, set $\s=\sin$ and $\c=\cos$.
If $f$ is indefinite, set $\s=\sinh$ and $\c=\cosh$.

\begin{prop}
\label{s4-param}
For a nondegenerate quadratic form $f=mx_1^2+kx_1x_2+nx_2^2$ such that
$mn\ne 0$, all possible integer
quadruples $(a,b,c,d)$ defining an integer normed pairing of type $(-,-)$ with
respect to $f$, belong to the
following family in $\R^4$ depending on one real parameter $\theta$:
\begin{equation}
\label{quad}
\begin{array}{l}
  a=\pm \sqrt{m} \frac{\ \s(3\theta+\varphi)}{\s(\varphi)} \\
  b=\pm \frac{n}{\sqrt{m}} \frac{\ \s(\theta + \varphi)(4 \c^2(\theta+\varphi)-1))}{\s (\varphi)} \\
  c =\pm \sqrt{n} \frac{\ \s(3\theta+ 2 \varphi)}{\s(\varphi)} \\
  d=\pm \frac{m}{\sqrt{n}} \frac{\ \s (\theta)(4 \c^2(\theta)-1))}{\s(\varphi)} \\
 \end{array} \end{equation}
Here $\varphi={\rm sign}(k)\ \c^{-1}(k/\sqrt{4mn})$.
All signs are equal: they are all pluses, or they are all minuses.
\end{prop}

\proof
Consider an orientation preserving real vector space isomorphism
$\phi:\R^2\to\A$ such that $f=|\phi|^2$.
Suppose that the matrix of $\phi$ with respect to the bases
$(\e_1,\e_2)$ of $\R^2$ and $(1,j)$ of $\A$ is
$$
\left(\begin{array}{cc}
\alpha&\beta\\
\gamma&\delta
 \end{array}\right).
$$
In coordinates, the equation $f=|\phi|^2$ reads
\begin{equation}
\label{in-L}
\alpha^2+\epsilon\gamma^2=m,\quad
\beta^2+\epsilon\delta^2=n,\quad
2(\alpha \beta+\epsilon\gamma\delta)=k,
\end{equation}

The general solution of this equation is
\begin{equation}
\label{abgd}
 \begin{array}{l}
  \alpha=\pm \sqrt{m}\ \c(\theta),\quad
  \beta=\pm\sqrt{n}\ \c(\theta+\varphi),\\
 \gamma=\pm \sqrt{m}\ \s(\theta),\quad
  \delta=\pm\sqrt{n}\ \s(\theta+\varphi),\\
 \end{array}
\end{equation}
where  $\varphi$ is defined by the equation $k=2\sqrt{mn}\ \c(\varphi)$.
The signs are the same: either they are all pluses or they are all minuses.
It is clear that any two orientation preserving isomorphisms $\phi$ with the property
$f=|\phi|^2$ differ by a rotation of $\A$ (see Theorem \ref{FormvsLat}).
In the case $\A=\H$, the choice of a sign corresponds to the choice of a
connected component in the rotation group $SO(1,1)$.

The numbers $(a,b,c,d)$ can be easily found in terms of $(\alpha,\beta,\gamma,\delta)$
by the formula
$$
s(\x,\y)=\phi^{-1}\sigma_4(\phi(\x),\phi(\y)),\quad \x,\y\in\R^2.
$$
Namely, we have
\begin{equation}
\label{abcd}
\begin{array}{ll}
 a=(\delta(\alpha^2-\epsilon \gamma^2)+2\alpha\beta\gamma)/(\alpha\delta-\beta\gamma),
 &
  b= \delta(3 \beta^2 -\epsilon \delta^2)/(\alpha \delta-\beta\gamma), \\
  c= (\gamma(\beta^2-\epsilon \delta^2)+2\alpha\beta\delta)/(\alpha \delta-\beta
  \gamma), &
  d= \gamma( 3\alpha^2-\epsilon \gamma^2)/(\alpha \delta-\beta \gamma).
 \end{array}
\end{equation}
The desired formula now follows.
$\Box$

\subsection*{Example 1}
Let us give an example of a nondegenerate integer quadratic form $f$ admitting
an integer normed pairings of types $(+,+)$, $(-,+)$ and $(-,-)$
but no integer normed pairing of type $(-,-)$.
Set $f=4x_1^2+2x_1x_2+6x_2^2$.
There are no values of $\theta$ such that all numbers
$a$, $b$, $c$, $d$ given by (\ref{quad}) are integer,
as Figure 1 shows.

\begin{figure}[ht]
%{\epsfbox{figinp.eps}}
%\epsffile{fig01.eps}
\includegraphics{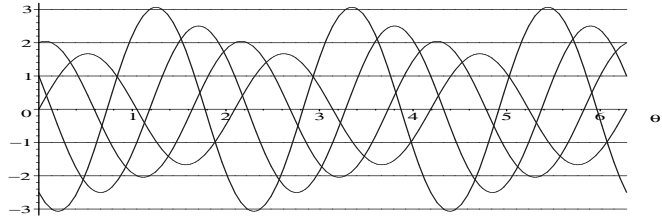}
\caption{Graphs of $a,c,b,d$: for $\theta=0$, $a=2,c=1,d=0,b=-5/2$}
\label{fig1}
\end{figure}

\subsection*{Example 2}
Let us give examples of nondegenerate integer quadratic forms $f$ admitting
integer normed pairings of all types $(+,+)$, $(-,+)$, $(+,-)$ and $(-,-)$.

We need the following fact, which is interesting on its own right:

\begin{thm}
\label{1-RDelta}
Suppose that an integer-normed lattice $L\subset\A$ contains 1.
Then $L=R_\Delta$, where $\Delta$ is the discriminant of $L$.
\end{thm}

\proof
Let $(1,\zeta)$ be a basis of $L$.
The existence of such basis follows from Proposition \ref{basis}.
The lattice $L$ is stable under multiplication.
Indeed, by Theorem \ref{3grinter}, we have $zw=zw\overline{1}\in L$
for all $z,w\in L$.

In particular, we have $\zeta^2\in L$.
Therefore, $\zeta^2$ is a linear combination of 1 and $\zeta$ with integer coefficients:
$$
\zeta^2+b\zeta+c=0.
$$

Solving this quadratic equation for $\zeta$, we obtain
$$
\zeta=\frac{-b+j\sqrt{|\Delta|}}2
$$
(we chose the plus sign to make the basis $(1,\zeta)$ properly oriented)
We see that $\Z[\zeta]=R_\Delta$.
But $L$ is stable under multiplication, hence $L=\Z[\zeta]$.
$\Box$

\begin{cor}
If an integer normed lattice $L\subset\A$ contains 1, then it is stable
under all operations $\sigma_k$, $k=1,\dots, 4$.
\end{cor}

\begin{cor}
If a nondegenerate integer quadratic form $f$ attains the value 1 at some
integer point, then $f$ admits integer normed pairings of all types
$(+,+)$, $(-,+)$, $(+,-)$ and $(-,-)$.
\end{cor}

This corollary provides many examples of quadratic forms admitting
integer normed pairings of all four types.
In particular, any form $x^2+Dy^2$, where $D$ is a nonzero integer,
admits integer normed pairings of all four types.

\subsection*{Example 3}
Let us give examples of nondegenerate integer quadratic forms $f$ admitting
integer normed pairings of type $(-,-)$ only.
To this end, we can use Theorem \ref{ord3}.

Recall that an integer quadratic form is called {\em primitive}, if
the greatest common divisor of all its coefficients is 1.
The following is an immediate corollary of Theorem \ref{inps}:

\begin{thm}
If a primitive integer quadratic form attains 1 nowhere on $\Z^2$,
then it admits no integer normed pairings of types $(+,+)$,
$(-,+)$ or $(+,-)$.
\end{thm}

Note that the primitivity condition is very essential here.
If $mp^2+kpq+nq^2\ne\pm 1$ in Theorem \ref{inps}, then
the forms $f_1$, $f_2$ and $f_3$ do not represent 1.
However, these forms are not primitive.

There are primitive integer quadratic forms attaining 1
nowhere on $\Z^2$ and admitting integer normed pairings
of type $(-,-)$.
Below are some particular examples of such forms, together with
their discriminants:
$$
\begin{array}{c}
f=(2,\pm 1,3),\quad \Delta=-23,\\
f=(2,\pm 1,4),\quad \Delta=-31,\\
f=(3,\pm 1,5),\quad \Delta=-59.
\end{array}
$$

\end{document}